\documentclass[leqno,12pt,letterpaper]{article}
\usepackage{exscale,here,latexsym,amsmath}
\usepackage{}
\usepackage{dsfont}
\long\def\ignore#1{}

\newcommand\comm[1]{\typeout{Used \string\comm...}\quad{\bf [[#1]]}\quad}

\newtheorem{theorem}{Theorem}

\newtheorem{lemma}[theorem]{Lemma}

\title{Dense arrangements are locally very dense I.}

\author{
J\'ozsef Solymosi\thanks{Department of Mathematics, University of
British Columbia, Vancouver, BC, Canada V6T 1Z2, Email:
\texttt{solymosi@math.ubc.ca} The research was supported by OTKA and
NSERC grants.}}

\begin{document}
\maketitle

\begin{abstract}
The Szemer\'edi-Trotter theorem \cite{st-epdg-83} gives a bound on
the maximum number of incidences between points and lines on the
Euclidean plane. In particular it says that $n$ lines and $n$ points
determine $O(n^{4/3})$ incidences. Let us suppose that an
arrangement of $n$ lines and $n$ points defines $cn^{4/3}$
incidences, for a given positive $c.$ It is widely believed that
such arrangements have special structure, but no results are known
in this direction. Here we show that for any natural number, $k,$
one can find $k$ points of the arrangement in general position such
that
 any pair of them is incident to a line from the arrangement, provided by $n\geq
 n_0(k).$ In a subsequent paper we will establish similar statement
 to hyperplanes.
\end{abstract}

\section{Introduction\label{sec:intro}}

The celebrated Szemer\'edi-Trotter Theorem \cite{st-epdg-83} states
that for $n$ points on the plane, the number of $m$-rich lines
cannot exceed
\begin{equation}
O(n^2/m^3+ n/m),
\end{equation}
and this bound is tight in the worst case. This result has numerous
applications not only in geometry \cite{ps-gi-04,s-cnhep-97}, but
also in number theory \cite{e-svp-02}. The Szemer\'edi-Trotter
theorem has various proofs, the most elegant is Sz\'ekely's
\cite{s-cnhep-97}. However the proofs provide very limited insight
view of the structure of extremal arrangements. It is widely
believed that a point-line arrangement which defines many incidences
has a special, somehow rigid structure. For example let us mention
here a question of Elekes. Is it true that for every $c>0$ there is
a $c'>0$ such that if a set of $n$ points on the plane contains at
least $cn^2$ collinear triples then at least $c'n$ points are along
an algebraic curve of degree $d,$ where $d$ is an universal
constant.

\medskip

The main purpose of this paper is to show that any arrangement with
close to the maximum number of incidences is locally a collection of
complete geometric graphs. For the sake of simplicity we state the
theorem for the balanced case, when the number of lines equals to
the number of points, but it is quite straightforward to see the
similar statement for unbalanced cases as well.

\medskip

Recent work of Gowers \cite{GO} and Nagle, R\"odl, Schacht, and
Skokan \cite{RO1,RO2,RO3} has established a hypergraph removal
lemma, which in turn implies similar results to hyperplanes, however
a slightly different approach is needed, mainly because the higher
dimensional extensions of the Szemer\'edi-Trotter theorem are not as
well defined as the planar case. To obtain sharp bounds one needs
certain restrictions on the arrangements. Therefore the
corresponding structure theorems will appear in subsequent paper.

\medskip

A pointset or a set of lines is in {\em general position} if no
three of the elements are collinear or concurrent.

\begin{theorem}
For every natural number $k$ and real $c>0$ there is a threshold
$n_0=n_0(k,c)$ such that if an arrangement of $n\geq n_0$ lines and
$n$ points defines at least $cn^{4/3}$ incidences, then one can
always find $k$ points of the arrangement in general position, such
that any pair of them is incident to a line from the arrangement.
\end{theorem}

As we will see it from the proof, the complete $k$-tuple is "local"
in a sense that for any pair of points, $p_1$ and $p_2,$ the number
of points from the arrangement, incident to the line segment
$(p_1,p_2)$ is less than $k.$

\section{Proof of Theorem 1}

The main tool of the proof is Szemer\'edi's Regularity Lemma
\cite{SZ1,SZ2}. We will use it's {\em counting lemma} form, because
it is easier to extend to hypergraphs which we will need for the
higher dimensional extensions. Let us prove first the simplest case,
to show that there is always a triangle. This "simplest case" is
interesting in it's own right, the statement of Lemma 2 implies
Roth's theorem \cite{RO} about arithmetic progressions on dense
subsets of integers. For the details we refer to \cite{SO1,SO2}.

\begin{lemma}        \label{lemma:line} For every $c>0$ there is a threshold
$n_0=n_0(c)$ and a positive $\delta =\delta (c)$ such that, for any
set of $n\geq n_0$ lines $L$ and any set of $m\geq cn^2$ points $P$,
if every point is incident to three lines, then there are at least
$\delta n^3$ triangles in the arrangement. (A triangle is a set of
three distinct points from $P$ such that any two are incident to a
line from $L.$) \end{lemma}

This lemma follows from the following theorem of Ruzsa and
Szemer\'edi \cite{RSZ}, which is also called as the {\em triangle
removal lemma} or the counting lemma for triangles.

\begin{theorem} \cite{RSZ} Let $G$ be a graph on $n$ vertices. If $G$ is the
union of $cn^2$ edge-disjoint triangles, then $G$ contains at least
$\delta n^3$ triangles, where $\delta$ depends on $c$ only.
\end{theorem}

The same theorem from a different angle is the following.

\begin{theorem}
Let $G$ be a graph on $n$ vertices. If $G$ contains $o(n^3)$
triangles, then one can remove $o(n^2)$ edges to make $G$
triangle-free.
\end{theorem}

To prove Lemma 2, let us construct a graph where $L$ is the vertex
set, and two vertices are adjacent if and only if the corresponding
lines cross at a point of $P$. This graph is the union of $cn^2$
disjoint triangles, every point of $P$ defines a unique triangle, so
we can apply Theorem 3.

\medskip

To determine the number of triangles in any arrangement of lines and
points seems to be a hard task. A related conjecture of de Caen and
Sz\'ekely \cite{CSZ} is that $n$ points and $m$ lines can not
determine more than $nm$ triangles.

\medskip

One can repeat the same argument, now with $k$ instead of 3. The
corresponding counting lemma can be proven using Szemer\'edi's
Regularity Lemma. The proof is analogous to the Ruzsa-Szemer\'edi
Theorem. There are slightly different ways to state the Regularity
Lemma, for our purposes the so called {\em degree form} is
convenient. For the notations and proofs we refer to the survey
paper of Koml\'os and Simonovits \cite{KS}.

\begin{theorem}{\bf (Regularity Lemma)}
For every $\epsilon > 0$ there is an $M = M(\epsilon )$ such that if
$G = (V,E)$ is any graph and $d \in (0,1]$ is any real number, then
there is a partition of the vertex-set $V$ into $k + 1$ clusters
$V_0, V_1,\ldots , V_k,$ and there is a subgraph $G'\subset G$ with
the following properties:

\begin{itemize}
\item $k \leq M,$
\item $|V_0|\leq \epsilon |V|,$
\item all clusters $V_i,$ $i\geq 1,$ are of the same size $m\leq\lceil \epsilon|V|\rceil ,$
\item $\deg_{G'}(v) > \deg_{G}(v) - (d+\epsilon)|V|$ for all $v\in
V,$
\item $e(G'(V_i))=0$ for each $i\geq 1,$
\item all pairs $G'(V_i, V_j)$ $(1\leq i < j \leq k)$ are $\epsilon-$regular, each with a density
either 0 or greater than $d$.
\end{itemize}
\end{theorem}

%\begin{theorem}{\bf (Counting Lemma for Graphs)}
%Let $G$ be a graph on $n$ vertices. If $G$ contains $o(n^k)$
%complete graphs on $k$ vertices, $K_k,$ then one can remove $o(n^2)$
%edges to make it $K_k-$free.
%\end{theorem}

Armed with the Regularity Lemma we are ready to prove the following
statement which is crucial in the proof of Theorem 1.

\begin{lemma}        \label{lemma:line} For every $c>0$ there is a threshold
$n_0=n_0(c)$ and a positive $\delta =\delta (c)$ such that, for any
set of $n\geq n_0$ lines $L$ and any set of $m\geq cn^2$ points $P$,
if every point is incident to $k$ lines, then there are at least
$\delta n^k$ complete $k-$tuples in the arrangement. (A complete
$k-$tuple is a set of $k$ distinct lines in general position from
$L$ such that any two intersect in a point from $P.$)
\end{lemma}

\begin{proof}
To avoid having too many degenerate $k-$tuples, we remove some
points from $P$ which have many lines incident to them. $P',$ which
is the subset of $P,$ consists of points incident to at most $100/c$
lines from $L.$ We can apply (1) to see that $P'$ is a large subset
of $P$, say $2|P'|>|P|.$ Let us construct a graph $G$ where $L$ is
the vertex set, and two vertices are adjacent if and only if the
corresponding lines cross at a point of $P'$. This graph, $G,$ is
the union of at least ${c\over 2}n^2$ edge-disjoint $K_k-$s. Find a
subgraph, $G',$ provided by Theorem 5 with $\epsilon \ll c.$ In $G'$
we still have some complete $K_k-$s. Along the edges of such a
complete graph we have a $k-$tuple of $V_i-$s such that the
bipartite graphs between them are dense and regular. This already
implies the existence of many complete subgraphs, as the following
lemma, quoted from \cite{KS}, shows.

\begin{lemma} Given $d > \epsilon > 0,$ a graph $R$ on $n$ vertices,
and a positive integer $m,$ let us construct a graph $G$ by
replacing every vertex of $R$ by $m$ vertices, and replacing the
edges of $R$ with $\epsilon-$regular pairs of density at least $d.$
Then $G$ has at least $\alpha m^n$ copies of $R,$ where $\alpha$
depends on $\epsilon, d,$ and $n,$ but not on $m.$
\end{lemma}

Most of the complete $k-$vertex subgraphs of graph $G'$ define a
complete $k-$tuple in the arrangement, i.e. the corresponding lines
are in general position. To see this, let us count the "degenerate"
$k-$tuples, the $k-$ element line-sets, where at least one triple is
concurrent. The number of concurrent triples is at most $cn^2{100/c
\choose 3}\leq c'n^2.$ For every concurrent triple one can select
$k-3$ lines to get a degenerate $k-$tuple. The expression,
$c'n^{k-1},$ is clearly an upper bound on the degenerate $k-$tuples,
therefore most of the complete graphs on $k$ vertices in $G'$ are
complete $k-$tuples if $n$ is large enough, $n\geq n_0=n_0(c).$

\end{proof}

The final step of the proof of Theorem 1 is to show that
arrangements with many incidences always have a substructure where
one use Lemma 2. We divide the arrangement into smaller parts where
we apply the dual of Lemma 2. The common technique to do that is the
so called {\em cutting}, which was introduced by Chazelle, see in
\cite{c-dm-00}, about 20 years ago. Here we use a more general
result, a theorem of Matousek \cite{m-ept-92}.

\begin{lemma} Let $P$ be a pointset, $P\subset \mathds{R}^d, |P|=n,$ and let
$r$ be a parameter, $1\ll r\ll n.$ Then the set $P$ can be
partitioned into $t$ sets $\Delta_1,\Delta_2,\ldots ,\Delta_t,$ in
such a way that $n/r\leq |\Delta_i|\leq 2n/r$ for all $i,$ and any
hyperplane crosses no more than $O(r^{1-1/d})$ sets, where $t=O(r).$
\end{lemma}

Here we use the $d=2$ case and we choose the value $r=\beta_k
n^{2/3},$ where $\beta_k$ is a constant, depends on $k,$ which we
will specify later.  Let us count the number of incidences along the
lines of $L,$ according to the partition of $P.$ For a given line
$\xi\in L,$ we count the sum $\sum_{i=1}^{t}\lfloor|\Delta_i\bigcap
\xi|/k\rfloor,$ which is not much smaller than the number of
incidences on $\xi$ over  $k$ if $\xi$ is rich of incidences, say,
incident to much more than $r^{1/2}k$ points of $P.$ From the
condition of Theorem 1 and the properties of the partition we have
the following inequality.
\begin{displaymath}
{c\over k}n^{4/3}\leq \sum_{\xi\in
L}\sum_{i=1}^{t}\left\lfloor{{|\Delta_i\bigcap
\xi|}\over{k}}\right\rfloor+|L|r^{1/2}
\end{displaymath}

Choosing $\beta_k={c\over 2k},$ the inequality becomes

\begin{displaymath}
{{cn^{4/3}}\over {2k}}=c_kn^{4\over3}\leq \sum_{\xi\in
L}\sum_{i=1}^{t}\left\lfloor{{|\Delta_i\bigcap
\xi|}\over{k}}\right\rfloor=\sum_{i=1}^{t}\sum_{\xi\in
L}\left\lfloor{{|\Delta_i\bigcap \xi|}\over{k}}\right\rfloor.
\end{displaymath}

Therefore there is an index $i$, such that

\begin{displaymath}
c_kn^{2/3}\leq\sum_{\xi\in L}\left\lfloor{{|\Delta_i\bigcap
\xi|}\over{k}}\right\rfloor.
\end{displaymath}

If $s=\left\lfloor{{|\Delta_i\bigcap \xi|}\over{k}}\right\rfloor$
then we can partition the points incident to $\xi$ into $r$
consecutive $k-$tuples. We can break the line into $r$ $k-$rich line
segments and consider them as separate lines. Our combinatorial
argument in Lemma 6 is robust enough to allow such modifications.
Then we have some $c'n^{2/3}$ $k-$rich lines on
$|\Delta_i|=c"n^{1/3}$ points. (Another possible way to show that
there are at least $c'n^{2/3}$ $k-$rich lines, is to apply the
Szemer\'edi-Trotter theorem, (1), to show that most of the lines are
not "very rich".) To complete the proof of Theorem 1, we apply the
dual statement of Lemma 6.


\begin{thebibliography}{10}

\topsep=0pt\itemsep=-2pt\parsep=1pt


\bibitem{KS} J.~Koml\'os, M.~Simonovits, Szemer\'edi's regularity lemma
and its applications in graph theory, in: Combinatorics, Paul Erd\H
os is eighty, Vol. 2 (Keszthely, 1993), 295–352, Bolyai Soc. Math.
Stud., 2, J\'anos Bolyai Math. Soc., Budapest, 1996.

\bibitem{RO1} B.~Nagle, V.~R\"odl, M.~Schacht, The counting lemma for
regular k-uniform hypergraphs, to appear, {\em Random Structures and
Algorithms.}

\bibitem{RO2} V.~R\"odl, J.~Skokan, Regularity lemma for k-uniform hypergraphs,
to appear, {\em Random Structures and Algorithms.}

\bibitem{RO3} V.~R\"odl, J.~Skokan, Applications of the regularity lemma for
uniform hypergraphs, preprint.

\bibitem{RO} K.F.~Roth, On certain sets of integers, {\em J. London Math. Soc.} 28
(1953), 245--252.

\bibitem{RSZ} I.~Ruzsa, E.~Szemer\'edi, Triple systems with no six points
carrying three triangles, {\em Colloq. Math. Soc. J. Bolyai} 18
(1978), 939–-945.

\bibitem{SO1} J.~Solymosi, Note on a generalization of Roth's theorem,
{in: Discrete and computational geometry}, 825–-827, Algorithms
Combin. 25, Springer Verlag, 2003.

\bibitem{SO2} J.~Solymosi, A note on a question of Erd\H os and Graham,
{\em Combinatorics, Probability and Computing} 13 (2004), 263-–267.

\bibitem{SZ1} E.~Szemer\'edi, On sets of integers containing no four elements
in arithmetic progression, {\em Acta Math. Acad. Sci. Hungar.} 20
(1969), 89–-104.

\bibitem{SZ2} E. Szemer\'edi, Regular partitions of graphs, in "Problem'es
Combinatoires et Th'eorie des Graphes, Proc. Colloque Inter. CNRS,"
(Bermond, Fournier, Las Vergnas, Sotteau, eds.), CNRS Paris, 1978,
399–-401.

\bibitem{c-dm-00} B.~Chazelle. The Discrepancy Method.
%: Randomness and Complexity.
Cambridge University Press, 2000.

\bibitem{ceg-ccbac-90}
K.~L.~Clarkson, H.~Edelsbrunner, L.~J.~Guibas, M.~Sharir, and
E.~Welzl, Combinatorial complexity bounds for arrangements of curves
and spheres, {\em Discrete Comput. Geom.} {\bf 5} (1990), 99--160.

\bibitem{e-svp-02}
Gy.~Elekes, Sums versus products in Number Theory, Algebra and Erd\H
os Geometry, vol 11 of Bolyai Soc. Math. Stud., J\'anos Bolyai Math.
Soc., Budapest, 2002, pp.~241--290.

\bibitem{et-inh-05}
Gy.~Elekes and Cs.~D.~T\'oth, Incidences of not-too-degenerate
hyperplanes, in {\em Proc. 21st ACM Sympos. Comput. Geom. (Pisa,
2005)}, ACM Press, to appear.


\bibitem{m-ept-92}
J.~Matou\v sek, Efficient partition trees, {\em Discrete Comput.
Geom.} {\bf 8} (1992), 315--334.

\bibitem{ps-gi-04}
J.~Pach and M.~Sharir, Geometric incidences, in {\em Towards a
Theory of Geometric Graphs}, vol.~342 of Contemporary Mathematics,
Amer. Math. Soc., Providence, RI, 2004, pp.~185--223.


\bibitem{st-ddp-01}
J.~Solymosi and Cs.~D.~T\'oth, Distinct distances in the plane, {\em
Discrete Comput. Geom.} {\bf 25} (2001), 629--634.



\bibitem{st-epdg-83}
E.~Szemer\'edi and W.~T.~T{rotter~Jr.}, Extremal problems in
{D}iscrete {G}eometry, {\em Combinatorica} {\bf 3} {\rm (3--4)}
(1983), 381--392.

\bibitem{s-cnhep-97}
L.~A.~Sz\'ekely, Crossing numbers and hard Erd\H os problems in
discrete geometry, {\em Combinatorics, Probability \& Computing}
{\bf 6} (3) (1997), 353--358.

\bibitem{GO} W.T.~Gowers, Hypergraph regularity and
the multidimensional Szemer\'edi theorem, preprint.

\bibitem{CSZ}
D.~de~Caen and L.~A.~Sz\'ekely, On dense bipartite graphs of girth
eight and upper bounds for certain configurations in planar
pointline systems, {\em J. Combin. Theory Ser. A.} {\bf 77} (1997),
268--278.

\end{thebibliography}
\end{document}